\documentclass[11pt]{article}

\usepackage{amsmath, amssymb, amscd}

\def\eqref#1{(\ref{#1})}
\newcommand{\goth}{\frak}

\newcommand{\arrow}{{\:\longrightarrow\:}}
\newcommand{\Z}{{\Bbb Z}}
\newcommand{\C}{{\Bbb C}}
\newcommand{\R}{{\Bbb R}}

\newcommand{\6}{\partial}
\def\1{\sqrt{-1}\:}
\newcommand{\restrict}[1]{{\left|_{{\phantom{|}\!\!}_{#1}}\right.}}

\newcommand{\calo}{{\cal O}}

\renewcommand{\bar}{\overline}
\renewcommand{\phi}{\varphi}
\renewcommand{\epsilon}{\varepsilon}
\renewcommand{\geq}{\geqslant}
\renewcommand{\leq}{\leqslant}

\newcommand{\const}{{\it const}}

\newcommand{\End}{\operatorname{End}}
\newcommand{\Tot}{\operatorname{Tot}}
\newcommand{\Id}{\operatorname{Id}}

\newcommand{\slope}{\operatorname{slope}}
\newcommand{\rk}{\operatorname{rk}}

\newcommand{\Tr}{\operatorname{Tr}}

\renewcommand{\Re}{\operatorname{Re}}

\newcommand{\comment}[1]{{}}

\def\blacksquare{\hbox{\vrule width 4pt height 4pt depth 0pt}}
\def\endproof{\blacksquare}

\makeatletter

\@ifundefined{Bbb}
     {\newcommand{\Bbb}[1]{{\mathbb #1}}}%
{}

\newcommand{\ps@verbit}{%
  \renewcommand{\@oddhead}{%
          \scriptsize
          {Stable bundles on elliptic fibrations}
          \hfil\tiny {M. Verbitsky, March 22, 2004 }}
  \renewcommand{\@evenhead}{\@oddhead}
  \renewcommand{\@oddfoot}{\hfil\thepage\hfil}
  \renewcommand{\@evenfoot}{\@oddfoot}}
 
\newcounter{Mycounter}[section]
\newcounter{lemma}[section]
\renewcommand{\thelemma}{\noindent{Lemma \thesection.\arabic{lemma}}}
\newcommand{\lemma}{%
     \setcounter{lemma}{\value{Mycounter}}
     \refstepcounter{lemma}
     \stepcounter{Mycounter}
     {\bf \thelemma:\ }}

\newcounter{claim}[section]
\newcounter{sublemma}[section]
\newcounter{corollary}[section]
\newcounter{theorem}[section]
\renewcommand{\thetheorem}{\noindent{Theorem \thesection.\arabic{theorem}}}
\newcommand{\theorem}{%
     \setcounter{theorem}{\value{Mycounter}}
     \refstepcounter{theorem}
     \stepcounter{Mycounter}
     {\bf \thetheorem:\ }}

\newcounter{conjecture}[section]
\newcounter{proposition}[section]
\renewcommand{\theproposition}
       {\noindent{Proposition \thesection.\arabic{proposition}}}
\newcommand{\proposition}{%
     \setcounter{proposition}{\value{Mycounter}}
     \refstepcounter{proposition}
     \stepcounter{Mycounter}
     {\bf \theproposition:\ }}

\newcounter{definition}[section]
\renewcommand{\thedefinition}
       {\noindent{Definition~\thesection.\arabic{definition}}}
\newcommand{\definition}{%
     \setcounter{definition}{\value{Mycounter}}
     \refstepcounter{definition}
     \stepcounter{Mycounter}
     {\bf \thedefinition:\ }}

\newcounter{example}[section]
\renewcommand{\theexample}{\noindent{Example \thesection.\arabic{example}}}
\newcommand{\example}{%
     \setcounter{example}{\value{Mycounter}}
     \refstepcounter{example}
     \stepcounter{Mycounter}
     {\bf \theexample:\ }}

\newcounter{remark}[section]
\renewcommand{\theremark}{\noindent{Remark \thesection.\arabic{remark}}}
\newcommand{\remark}{%
     \setcounter{remark}{\value{Mycounter}}
     \refstepcounter{remark}
     \stepcounter{Mycounter}
     {\bf \theremark:\ }}

\newcounter{problem}[section]
\newcounter{question}[section]
\@addtoreset{equation}{section}
\@addtoreset{footnote}{section}
\makeatother

\begin{document}

\begin{center}
{\LARGE\bf
Stable bundles on positive principal \\[3mm] elliptic fibrations 
}
\\[4mm]
Misha Verbitsky\footnote{Misha Verbitsky is 
an EPSRC advanced fellow 
supported by CRDF grant RM1-2354-MO02 
and EPSRC grant  GR/R77773/01}
\\[4mm]

{\tt verbit@maths.gla.ac.uk, \ \  verbit@mccme.ru}
\end{center}

{\small 
\hspace{0.15\linewidth}
\begin{minipage}[t]{0.7\linewidth}
{\bf Abstract} \\
Let $í\stackrel\pi \arrow X$ be a principal elliptic fibration 
over a K\"ahler base $X$. We assume that the
K\"ahler form on $X$ is lifted to an exact form
on $M$ (such fibrations are called {\bf positive}).
Examples of these are regular Vaisman
manifolds (in particular, the regular Hopf 
manifolds) and Calabi-Eckmann manifolds. 
Assume that $\dim M > 2$.
Using the Kobayashi-Hitchin correspondence, 
we prove that all stable bundles on $M$ are
flat on the fibers of the elliptic fibration.
This is used to show that all stable vector 
bundles on $M$ take form $L\otimes \pi^* B_0$,
where $B_0$ is a stable bundle on $X$, and $L$
a holomorphic line bundle. For $X$ algebraic 
this implies that
all holomorphic bundles on $M$ are filtrable
(that is, obtained by successive extensions
of rank-1 sheaves). We also show that
all positive-dimensional compact subvarieties
of $M$ are pullbacks of complex 
subvarieties on $X$.
\end{minipage}
}

{
\small
\tableofcontents
}


\section{Introduction}
\label{_Intro_Section_}


Let $M$ be a compact complex manifold, and $M \stackrel\pi\arrow X$
a smooth holomorphic fibration. Assume that an elliptic
curve (considered as a complex Lie group) 
acts on $M$ holomorphically. Assume, moreover,
that this action is free and transitive on the
fibers of $\pi$. Then $\pi:\; M \arrow X$ is called 
{\bf a principal elliptic fibration}. 
The fibers of $\pi$ are identified 
(non-canonically) with $T$.

This terminology is somewhat misleading, as in
algebraic geometry one says ``elliptic fibraion''
speaking of a fibration with a section. However,
a principal bundle with a section is trivial,
hence we use (following \cite{_Brinzanescu_Moraru:FM_},
\cite{_Brinzanescu_Moraru:stable_}) 
the term ``principal elliptic fibration'' used
to describe something which should be more properly 
called ``a 1-dimensional torus principal bundle''.

In physics, the same object is usually called ``a $T^2$-bundle''.
There is a great body of literature studying these manifolds
in physical context, see e.g. \cite{_Cardoso_etc_},
\cite{_Goldstein_Prokushin_} and the references therein.
One usually considers a principal elliptic fibraion
over a 2-dimensional Calabi-Yau manifold, equipped with
a Hermitian metric, such that its Bismut
connection has holonomy $SU(3)$.\footnote{Bismut
connection on $M$ is a Hermitian connection in $TM$,
with totally antisymmetric torsion; it is not difficult 
to show that such connection is unique on any
Hermitian complex manifold.} Such manifolds
appear as a background for the heterotic 
string in the presence of fluxes.

Examples of principal elliptic fibrations
are very common in complex geometry
(take, for instance, the Hopf surface
$\C^2\backslash 0 /\langle q \rangle$,
fibered over $\C P^1$ with a fiber 
$\C^*/\langle q \rangle$, where 
$q\in \C$, $|q|>1$). For more examples
see Subsection \ref{_exa_figra_Subsection_}.
The Hopf surface $M$ is clearly homeomorphic
to $S^1 \times S^3$. Since $H^2(M)=0$,
it is not K\"ahler. This is a general 
phenomenon. Unless $\pi$ is 
trivialized after taking a finite 
cover, the total space $M$ of a
principal elliptic fibration
is not K\"ahler. 

In this paper we are interested in the so-called
{\bf positive principal elliptic fibrations} 
(\ref{_positive_fibra_Definition_}).
These are the fibrations with the K\"ahler base $X$,
such that some K\"ahler form $\omega_X$ on $X$ is lifted
to an exact form on $M$. Heuristically this means
that the fibration $M \stackrel\pi\arrow X$
is positive in the sense of complex algebraic
geometry (see \ref{_posi_for_Vaisman_Remark_}  
for a more detailed explanation). From the positivity
one obtains that $M$ admits an exact positive
2-form $\pi^*\omega_X$ with all eigenvalues
strictly positive except one. 

In \cite{_Verbitsky:LCHK_}, \cite{_OV:Immersion_}
such form was used in Vaisman geometry
to obtain Ko\-da\-ira-\-type embedding theorem,
several Kodaira-Nakano-type vanishing results and
a study of compact subvarieties of a Vaisman manifold. We
apply a similar reasoning to positive elliptic fibrations.
Using the Kobayashi-Hitchin correspondence 
(Subsection \ref{_Koba_Hi_Subsection_}), we obtain
the following theorem.

\hfill

\theorem\label{_curva_vani_Intro_Theorem_}
Let $M\stackrel \pi\arrow X$, $\dim_\C M =n$ be a positive principal
elliptic fibration, $(B, \nabla)$ a Hermitian-Einstein bundle on $M$, and
$\Theta\in \Lambda^{1,1}(M)\otimes \End(B)$ its curvature.
Assume that $n\geq 3 $. Then $\Theta(v, \cdot)=0$
for any vertival tangent vector $v\in T_\pi M$.

\hfill

{\bf Proof:} This is \ref{_curva_YM_vanishes_verti_Theorem_}.
\endproof

\hfill

For a definition of stability and Hermitian-Einstein
connections see Section \ref{_stabi_Koba_Hi_Section_}.
\ref{_curva_vani_Intro_Theorem_} implies the following
corollary.

\hfill

\proposition
Let $T$ be an elliptic curve, and 
$M \stackrel \pi \arrow X$ a positive principal
$T$-fibration, equipped with a preferred Hermitian
metric. The universal covering $\tilde T$ acts
on $M$ in a standard way (its action is factorized 
through $T$). Consider a stable bundle $B$
on $M$. Then $B$ is equipped with a natural
holomorphic $\tilde T$-equivariant structure.

\hfill

{\bf Proof:} This is \ref{_equiva_sta_Proposition_}.
\endproof

\hfill

A similar argument proves

\hfill

\proposition
Let $M\stackrel \pi\arrow X$ be a positive principal
elliptic fibration, and $Z\subset M$ a closed positive-dimensional
subvariety. Then $Z= \pi^*(Z_0)$ for some complex subvariety
$Z_0\subset X$.

\hfill

{\bf Proof:} This is \ref{_subva_Proposition_}. \endproof

\hfill

Applying the $\tilde T$-equivariant structure arising on stable
bundles, the following structure theorem is obtained

\hfill

\theorem
Let $M\stackrel\pi\arrow X$ be a positive principal
elliptic fibration equipped with a preferred Hermitian
metric, and $B$ a stable holomorphic bundle
on $M$. Then $B\cong L\otimes \pi^* B_0$, where
$L$ is a line bundle on $M$ and $B_0$ a stable
bundle on $X$.

\hfill

{\bf Proof:} This is \ref{_stru_sta_Theorem_}. \endproof

\hfill

We also show that all coherent sheaves on $M$ are filtrable, that
is, obtained as successive extensions of rank 1 sheaves
(\ref{_filtrable_Theorem_}).


\section{Principal elliptic fibrations}


\subsection{Positive elliptic fibrations: definition and examples}
\label{_exa_figra_Subsection_}

Throughout this paper, $M$ is a compact complex manifold,
$X$ a K\"ahler manifold, $T$ an elliptic curve, and
$M\stackrel \pi \arrow X$ a principal $T$-fibration

\hfill

\definition\label{_positive_fibra_Definition_}
A fibration is called {\bf positive} if the pullback
$\pi^*\omega_X$ is exact, for some K\"ahler 
form $\omega_X$ on $X$.

\hfill

Principal toric fibrations, their invariants, 
topology and Dolbeault cohomology are 
thoroughly analyzed in the excellent paper 
\cite{_Hofer:remarks_}. For our purposes, the
most important examples are the following.

\hfill

\example\label{_regu_Vaisman_via_ample_Example_}
Regular Vaisman manifolds are principal
elliptic fibrations constructed the following
way (see e.g. \cite{_Dragomir_Ornea_}, 
\cite{_OV:Immersion_}). Take a projective
manifold $X$, and let $L$ be an ample
line bundle on $X$. Consider the space
$\tilde M := \Tot(L^* \backslash 0)$,
which is a total space of the dual bundle
$L^*$ without the zero section. Then $\tilde M$
is a principal $\C^*$-bundle over $X$.
Fix $q\in \C$, $|q|>1$, and let $M:= \tilde M /\sim_q$,
be the quotient of $\tilde M$ under the equivalence
$v \sim q v$, $v \in L^*\backslash 0$.
Then $M$ is a principal elliptic bundle,
with a fiber $T = \C^* /\langle q \rangle$. 
The positivity of $M$ is elementary
(see e.g. \cite{_Verbitsky:LCHK_} or
\cite{_OV:Immersion_}).

\hfill

Vaisman manifolds were studied in great detail
by I. Vaisman (see e.g. \cite{_Vaisman:Dedicata_}),
under the name of generalized Hopf manifolds. 

\hfill

\remark 
A special case of the above example is a regular
Hopf manifold $M \cong S^{2n-1}\times S^1$, obtained as a quotient of
$\C^n\backslash 0$ under an equivalence
generated by $v \arrow qv$. Clearly,
$M$ is fibered over $\C P^{n-1}$,
with a fiber $T = \C^* /\langle q \rangle$.
Regular Hopf manifold 
is obtained  if one applies the construction of
\ref{_regu_Vaisman_via_ample_Example_}
to $X = \C P^{n-1}, L = \calo(1)$.

\hfill

\remark\label{_posi_for_Vaisman_Remark_} 
Taking an arbitrary line bundle $L$ on a K\"ahler manifold
$X$ and applying the same construction as above, we 
obtain a principal toric bundle as well. It will be
positive in the sense of \ref{_positive_fibra_Definition_}
if and only if the bundle $L$ is positive, as one
can easily see from the exact sequence associated
with the fibration in \cite{_Hofer:remarks_} 
\begin{equation}\label{_fibra_cohomo_exa_Equation_}
H^1(T) \stackrel \delta \arrow H^2(X) \stackrel i \arrow H^2(M).
\end{equation}
The second arrow of \eqref{_fibra_cohomo_exa_Equation_}
is a standard pullback map, and the transgression 
$\delta$ maps one of the generators of $H^1(T)$ 
to zero, the other to $c_1(L)$. 

\hfill

\example
The Calabi-Eckmann manifolds $M \cong S^{2n+1}\times S^{2m+1}$
are fibered in standard way over $\C P^n \times \C P^m$
(the spheres $S^{2n+1}$, $S^{2m+1}$ are fibered 
over $\C P^n$, $\C P^m$ by the way of the Hopf fibration). 
The complex structure on $M$ is constructed as follows. 
Fix a number $\tau \in \C\backslash \R$. 
Consider the action of an additive group $\C$ on 
$\C^{n+1} \backslash 0)\times (\C^{m+1} \backslash 0$,
\[ 
   t(v_1, v_2) \arrow (e^t  v_1, e^{t\tau}  v_2).
\]
This action is clearly holomorphic, and the quotient
space $M$ is naturally identified with 
$S^{2n+1}\times S^{2m+1}$;\footnote{To see this isomorphism,
take a Hermitian metric on $\C^{n+1}$, $\C^{m+1}$.
Then $|e^t v_1| = |e^{t\tau}v_2| =1$ if and only if
$\Re(t) = \log|v_1|, \Re(\tau t) = \log|v_2|$.
Since $\tau \notin \R$, such $t$ exists and is clearly unique.}
moreover, $M$ is equipped with a natural holomorphic
projection  $M \stackrel \pi \arrow \C P^n \times \C P^m$. 
The fibers of $\pi$ are identified with 
$\C/\langle 1, \tau\rangle$. 
It is not difficult to see that $M$
is a principal elliptic fibration.
Since $H^2(M)=0$, $M$ is obviously
positive. 

\subsection{Preferred Hermitian metrics}

Let  $M \stackrel \pi \arrow X$ be a  principal
elliptic fibration. As always, we assume that
$M$ is compact. By Blanchard's Theorem, unless
$\pi$ is trivialized on a finite covering, $M$ is
non-K\"ahler (see e.g. \cite{_Hofer:remarks_}). 
If $M$ is positive, the pullback of a K\"ahler 
metric becomes exact, hence $\pi$ cannot
be trivialized. Therefore, a positive principal
elliptic fibration is never K\"ahler.

We shall always assume that $M$ is equipped with a special
kind of Hermitian metric, defined as follows. 

\hfill

Recall that the smooth surjective 
map $\pi:\; (M, g) \arrow (X, g_X)$ of 
Riemannian manifolds is called {\bf Riemannian submersion}
if its differential $d\pi$ induces an isometry 
$d\pi:\; T_\pi M^\bot \arrow TX$,
where $T_\pi M$ denotes the space of vertical
tangent vectors, and $T_\pi M^\bot$ 
is orthogonal complement.

\hfill

\definition 
Let $T$ be an elliptic curve, and
$M \stackrel \pi \arrow X$ a positive principal
$T$-fibration. Consider a K\"ahler metric
$g_X$ on $X$, such that the pullback $\pi^* \omega_X$
of the corresponding K\"ahler form is exact. 
A Hermitian metric $g$ on $M$ is called 
{\bf preferred} if $g$ is $T$-invariant, and the projection
$\pi:\; (M, g) \arrow (X, g_X)$ is a Riemannian submersion.

\hfill

Clearly, preferred Hermitian metrics always exist.

\hfill

\remark
Let $M$ be a regular Vaisman manifold from 
\ref{_regu_Vaisman_via_ample_Example_}, and $\tau\in TM$ the real vector
field which is tangent to the map $v \arrow \lambda v$, 
$\lambda \in \R$, $v\in L^*$. Then $M$ admits a 
Hermitian metric $g$ for which $\tau$ acts by isometries,
and, moreover, $\nabla_{LC} \tau =0$ for the Levi-Civita
connection associated with this metric. In terminology
used in Vaisman geometry, $g$ is called the
{\bf locally conformally K\"ahler Gauduchon metric},
and $\tau$ {\bf the Lee field} (\cite{_Dragomir_Ornea_}).
From the condition $\nabla_{LC} \tau =0$ one immediately obtains
that $\pi:\; (M, g) \arrow (X, g_X)$ is a Riemannian
submersion (see e.g. \cite{_OV:Structure_}). 
By construction, $\tau$ is a generator of the Lie algebra
of $T$. Since $\tau$ acts by isometries, $g$ is 
$T$-invariant.  We obtain that the locally 
conformally K\"ahler Gauduchon metric on a 
regular Vaisman manifold is a preferred one.


\section{Stable bundle on Hermitian
manifolds}
\label{_stabi_Koba_Hi_Section_}


In this section, we recall the necessary definitions
and explain the Kobayashi-Hitchin correspondence on
non-K\"ahler manifolds. We follow \cite{_Lubke_Teleman:Book_}
and \cite{_Lubke_Teleman:Universal_}.

\subsection{Gauduchon metrics and stability}

\definition 
Let $(M, g)$ be a Hermitian manifold, and 
$\omega\in \Lambda^{1,1}(M)$ its Hermitian form.
We say that $g$ is a Gauduchon metric if
$\6^* \bar\6^* \omega=0$, or, equivalently,
$\6\bar\6(\omega^{\dim_\C M-1})=0$.

\hfill

Recall that the metrics $g$, $g'$ on $M$ are called {\bf conformally
equivalent} if $g= f g'$, $f\in C^\infty M$.

\hfill

\theorem 
(\cite{_Gauduchon_1984_})
Let $(M,g)$ be a compact Hermitian manifold.
Then there exists a unique Gauduchon metric $g'$ which is
conformally equivalent to $g$.

\endproof

\hfill

\definition
Let $M$ be a compact complex manifold equipped 
with a Gauduchon metric, and $\omega\in \Lambda^{1,1}(M)$ 
the corresponding Hermitian form. Consider a torsion-free
coherent sheaf $F$ on $M$. Denote by $\det F$ its determinant
bundle. Pick a Hermitian metric $\nu$ on $\det F$, and
let $\Theta$ be the curvature of the associated Chern
connection. We define the degree of $F$ as follows:
\[
\deg F := \int_M \Theta \wedge \omega^{\dim_\C M-1}.
\]
This notion is independent from the choice of the metric
$\nu$. Indeed, if $\nu' = e^\psi \nu$, $\psi \in C^\infty(M)$,
then the associated curvature form is written as
$\Theta' = \Theta + \6\bar\6 \psi$, and
\[ \int_M \6\bar\6 \psi\wedge \omega^{\dim_\C M-1}=0
\]
because $\omega$ is Gauduchon.

\hfill

If $F$ is a Hermitian vector bundle, $\Theta_F$ its
curvature, and the metric $\nu$ is induced from $F$, then
$\Theta= \Tr_F\Theta_F$. In K\"ahler case this allows
one to relate degree of a bundle with the first Chern class.
However, in non-K\"ahler case, the degree is not a 
topological invariant --- it depends fundamentally 
on the holomorphic geometry of $F$. 
Moreover, the degree is not discrete, as in the K\"ahler 
situation, but takes values in continuum.

Further on, we shall
see that one can in some cases construct 
a holomorphic structure of any given 
degree $\lambda\in \R$ on a fixed 
$C^\infty$-bundle. In our examples, such
holomorphic structures are constructed
on a topologically trivial line bundle 
over a positive principal elliptic 
fibration (\ref{_degree_triv_Remark_}).

\hfill

\definition
Let $F$ be a non-zero torsion-free coherent sheaf on $M$.
Then $\slope(F)$ is defined as 
\[
\slope(F) := \frac{\deg F}{\rk F}.
\]
The sheaf $F$ is called\\
{\small
\begin{tabular}[t]{ll}
{\bf stable} & if for all subsheaves $F'\subset F$, 
we have $\slope(F')< \slope(F)$\\
{\bf semistable} & if for all subsheaves $F'\subset F$, 
we have $\slope(F')\leq \slope(F)$\\
{\bf polystable} & if $F$ 
can be represented as a direct sum of stable \\& coherent 
sheaves with the same slope. 
\end{tabular}
}

\hfill

\remark
This definition is stability is ``good'' as most
standard properties of stable and semistable bundles 
hold in this situation as well. In particular, all line bundles 
are stable; all stable sheaves are simple; the Jordan-H\"older
and Harder-Narasimhan filtrations are well defined
and behave in the same way as they do in
the usual K\"ahler situation 
(\cite{_Lubke_Teleman:Book_}, \cite{_Bruasse:Harder_Nara_}).

However, not all bundles are {\bf filtrable},
that is, are obtained as successive extensions
by coherent sheaves of rank 1. There are non-filtrable
vector bundles e.g. on a non-algebraic K3 surface 
(\ref{_non_filtra_Example_}). 

\subsection{Kobayashi-Hitchin correspondence}
\label{_Koba_Hi_Subsection_}

The statement of Kobayashi-Hitchin correspondence is translated
to the Hermitian situation verbatim, following 
Li and Yau (\cite{_Li_Yau_}).

\hfill

\definition
Let $B$ be a holomorphic Hermitian vector
bundle on a Hermitian manifold $M$, and 
$\Theta\in \Lambda^{1,1}(M)\otimes \End(B)$
the curvature of its Chern connection $\nabla$. Consider the
operator $\Lambda:\; \Lambda^{1,1}(M)\otimes \End(B)\arrow \End(B)$
which is a Hermitian adjoint to $b \arrow \omega\otimes b$,
$\omega$ being the Hermitian form on $M$.
The connection $\nabla$ is called {\bf Hermitian-Einstein}
(or {\bf Yang-Mills}) if $\Lambda \Theta = \const \cdot \Id_B$.

\hfill

\theorem 
(Kobayashi-Hitchin correspondence)
Let $B$ be a holomorphic vector bundle on a compact 
complex manifold equipped with a Gauduchon metric. Then
$B$ admits a Hermitian-Einstein connection $\nabla$
if and only if $B$ is polystable. Moreover, the
Hermitian-Einstein connection is unique.

\hfill

{\bf Proof:} See \cite{_Li_Yau_}, \cite{_Lubke_Teleman:Book_}, 
\cite{_Lubke_Teleman:Universal_}. \endproof

\hfill

\remark 
The bundle $B$ is stable
if and only if the corresponding Hermitian-Einstein
bundle $(B, \nabla)$ cannot be decomposed onto a 
direct sum of sub-bundles. In this case,
the Hermitian-Einstein metric is defined by
the connection uniquely up to a constant
multiplier.


\section{Hermitian-Einstein bundles on 
positive principal elliptic fibrations}


\subsection{Preferred metrics are Gauduchon}

\proposition
Let $M\stackrel \pi\arrow X$ be a positive principal
elliptic fibration, and $g$ a preferred Hermitian
metric. Then $g$ is Gauduchon.

\hfill

{\bf Proof:} Consider the orthogonal decomposition
$TM = T_\pi M \oplus T_\pi M^\bot$
of the tangent bundle onto its vertical
component and its orthogonal complement.
We write the Hermitian form $\omega$ as
$\omega=\omega_0+\omega_T$, where
$\omega_0$ vanishes on $T_\pi M$ and
$\omega_T$ vanishes on $T_\pi M^\bot$. Since $\pi$ is a 
Riemannian submersion, $\omega_0=\pi^*\omega_X$, where
$\omega_X$ is a K\"ahler form on $X$.
This implies 
\[ \6^*\bar\6^*\omega_0 = \pi^*\6^*\bar\6^*\omega_X=0.\]

On the other hand, $\6^*\bar\6^*\omega_T$
can be computed by restricting $\omega_T$ to the
fibers $\pi^{-1}(x)$ of $\pi$ and computing 
$\6^*\bar\6^*(\omega_T\restrict{\pi^{-1}(x)})$.
Since $\omega_T$ is $T$-invariant, 
$\6^*\bar\6^*(\omega_T\restrict{\pi^{-1}(x)})$
also vanishes, and we have
\[ 
    \6^*\bar\6^*\omega = \6^*\bar\6^*\omega_0+\6^*\bar\6^*\omega_T=0.
\]
\endproof

\subsection{Primitive forms on positive principal elliptic
fibrations}

Let $B_0$ be a vector bundle with connection,
and $\eta$ a $B_0$-valued form. Then $\eta$ is called 
{\bf closed} if $\nabla\eta=0$, where 
\[ \nabla:\; \Lambda^i(M) \otimes B \arrow 
   \Lambda^{i+1}(M) \otimes B
\]
is the connection operator extended to forms
by the Leibniz rule. The curvature form of a connection
is closed, by Bianchi identity.

Further on in this section, we shall need the 
following result.

\hfill

\proposition\label{_primi_form_Proposition_}
Let $M\stackrel \pi \arrow X$ be a positive principal elliptic
fibration, $n = \dim M \geq 3$, 
equipped with a preferred Hermitian metric,
$B$ a Hermitian bundle with connection, and
$\Theta \in \Lambda^{1,1}(M) \otimes \End(B)$
a closed (1,1)-form. Assume that $\Theta$ 
is primitive, that is, $\Lambda\Theta=0$.
Then $\Theta(v, \cdot) =0$ for any vertcal
tangent vector $v\in T_\pi M$.

\hfill

{\bf Proof:} Since the connection $\nabla$ is Hermitian,
it preserves the natural real structure in 
$\Lambda^{1,1}(M) \otimes \End(B)$, 
$\eta \otimes b \arrow \bar\eta \otimes b^\bot$,
where by $b^\bot$ one understands the Hermitian 
adjoint endomorphism. Therefore, we may assume that
$\Theta$ is real, with respect to this real structure.

Let $\theta$, $\theta_1, ... , \theta_{n-1}$ be an orthonormal
basis in $\Lambda^{1,0}(M)$, with $\theta\in T_\pi M$, 
$\theta_i\in T_\pi M^\bot$. Consider a decomposition
\begin{align}\label{_Theta_basis_Equation_}
\Theta &= \sum_{i\neq j}(\theta_i \wedge \bar\theta_j
         + \bar\theta_i \wedge \bar\theta_j) \otimes b_{ij} +
\sum_{i}(\theta_i \wedge \bar\theta_i) \otimes a_i \\
& + \sum_{i}(\theta \wedge \bar\theta_i
         + \bar\theta \wedge \bar\theta_i) \otimes b_{i}
  + \theta\wedge\bar\theta\otimes a,
\end{align}
with $b_{ij}$, $b_i$, $a_i$, $a\in {\goth u}(B)$
being skew-Hermitian endomorphisms of $B$. 

Consider now the form $\omega_0:= \pi^* \omega_X$.
This form is exact, positive, and has $n-1$
strictly positive eigenvalues. Using the basis described above,
we can write
\[
\omega= \1 
\left(\theta\wedge\bar\theta +\sum_{i}\theta_i \wedge \bar\theta_i \right), 
\ \ \  \omega_0 = \1 \left (\sum_{i}\theta_i \wedge \bar\theta_i\right )
\]
where $\omega$ is the Hermitian form of $M$.
This follows directly from $\omega$ being a preferred
Hermitian form. 

Let $\Xi:= \Tr(\Theta\wedge \Theta)$. 
This is a closed (2,2)-form on $M$.
Then \eqref{_Theta_basis_Equation_} implies
\[
 (\1)^n\Xi \wedge \omega_0^{n-2} = \Tr\left(-\sum b_i^2 +a \left(\sum a_i\right)\right)
\]
On the other hand, $ \sum a_i + a = \Lambda\Theta=0$,
hence 
\[
 (\1)^n\Xi \wedge \omega_0^{n-2} = \Tr\left(-\sum b_i^2 -a^2 \right).
\]
Since $\Tr(-a^2)$ is a positive definite form on
${\goth u}(B)$, the integral
\begin{equation}\label{_integral_Xi_Equation_}
\int_M  (\1)^n\Xi \wedge \omega_0^{n-2}
\end{equation}
is non-negative, and positive unless $b_i$ and $a$ 
both vanish everywhere. Using $n>2$, we find that
\eqref{_integral_Xi_Equation_} vanishes, because
$\omega_0$ is exact and $\Xi$ is closed. Therefore,
$b_i$ and $a$ are identically zero, which is exactly the
claim of \ref{_primi_form_Proposition_}
\endproof

\subsection{The L\"ubke-type positivity for Hermitian-Einstein bundles}

\theorem\label{_curva_YM_vanishes_verti_Theorem_}
Let $M\stackrel \pi\arrow X$, $\dim_\C M =n$ be a positive principal
elliptic fibration, equipped with a preferred Hermitian metric,
$(B, \nabla)$ a Hermitian-Einstein bundle on $M$, and
$\Theta\in \Lambda^{1,1}(M)\otimes \End(B)$ its curvature.
Assume that $n\geq 3$. Then $\Theta(v, \cdot)=0$
for any vertical tangent vector $v\in T_\pi M$.

\hfill

{\bf Proof:} \ref{_curva_YM_vanishes_verti_Theorem_} is proven
by the standard positivity argument, going back to M. L\"ubke
(\cite{_Lubke_}). Similar argument is used e.g. to show 
that any Hermitian-Einstein bundle with vanishing Chern 
classes is flat.

Let $\Theta_0:= \Theta - \frac 1{\rk F}\Tr \Theta$
be the traceless part of $\Theta$. Then $\Theta_0$ is primitive
and closed, hence, by \ref{_primi_form_Proposition_},
$\Theta_0(v, \cdot)=0$ for all
$v \in T_\pi M$.
To prove \ref{_curva_YM_vanishes_verti_Theorem_},
it remains to show that the 2-form $\Tr \Theta$
also vanishes on $v \in T_\pi M$.
This is implied by the following trivial 
lemma.

\hfill

\lemma\label{_closed_forms_Lambda_const_vanish_on_T_pi_Lemma_}
Let $M\stackrel \pi\arrow X$, $\dim_\C M \geq 3$ be a positive principal
elliptic fibration, equipped with a preferred Hermitian metric,
and $\eta\in \Lambda^{1,1}(M)$ a closed form satisfying 
$\Lambda \eta= \const$. Then $\eta(v, \cdot)=0$ for all
$v \in T_\pi M$.

\hfill

{\bf Proof:} 
Clearly, $\Lambda \omega_0=\const$.
Replacing $\eta$ with $\eta-c \omega_0$, we may
assume that $\Lambda \eta=0$. Now $\eta$ is primitive,
and, by \ref{_primi_form_Proposition_}, 
$\eta(v, \cdot)=0$ for all $v \in T_\pi M$.
This finishes the proof of
\ref{_curva_YM_vanishes_verti_Theorem_}.
\endproof

\subsection{Complex subvarieties in positive principal elliptic fibrations}

A similar argument proves the following result.

\hfill

\proposition\label{_subva_Proposition_}
Let $M\stackrel \pi\arrow X$ be a positive principal
elliptic fibration, and $Z\subset M$ a closed positive-dimensional
subvariety. Then $Z= \pi^*(Z_0)$ for some complex subvariety
$Z_0\subset X$.

\hfill

{\bf Proof:} The proof is taken from \cite{_Verbitsky:LCHK_}
verbatim, where the same result was proven for
Vaisman manifolds.
 
 It is well known that
\begin{equation}\label{_int_of_omega_0_Equation_}
\int_Z \omega_0^k\geq 0
\end{equation}
for all complex subvarieties
$Z\subset M$, $\dim_\C Z =k$, and all positive forms $\omega_0$. 
Moreover, the integral \eqref{_int_of_omega_0_Equation_}
vanishes if and only if $Z$ is tangent to the
null-space foliation of $\omega_0$. 

Since $\omega_0$ is exact, the integral 
\eqref{_int_of_omega_0_Equation_}
vanishes. Therefore, $Z$ is tangent to the
null-space foliation of $\omega_0$. This implies
$Z$ is tangent to the fibers of $\pi$.
\endproof

\subsection{Equivariance of stable bundles}

\ref{_curva_YM_vanishes_verti_Theorem_} has an immediate
implication for the holomorphic geometry of stable bundles.

\hfill

\proposition\label{_equiva_sta_Proposition_}
Let $T$ be an elliptic curve, and 
$M \stackrel \pi \arrow X$, $\dim_\C M \geq 3$
a positive principal
$T$-fibration, equipped with a preferred Hermitian
metric. The universal covering $\tilde T$ acts
on $M$ in a standard way (its action is factorized 
through $T$). Consider a stable bundle $B$
on $M$. Then $B$ is equipped with a natural
holomorphic $\tilde T$-equivariant structure.

\hfill

{\bf Proof:} Let $\nabla$ be a Hermitian-Einstein
connection on $B$, and $\tau$ a holomorphic vector field
tangent to the action of $T$. Consider the operator
$\nabla_\tau:\; B \arrow B$. Since the curvature of
$B$ vanishes on $\tau$, $\nabla_\tau$ is holomorphic.
For any $\lambda\in \C$
the operator $e^{\lambda \nabla_\tau}$ of parallel
translation along $\tau$ defines a
holomorphic $\tilde T$-equivariant structure
on $B$. \endproof

\hfill

\remark
The same argument implies that a stable bundle on a Vaisman
manifold is always equivariant with respect to the action
of the Lee field (see e.g. \cite{_Dragomir_Ornea_},
\cite{_OV:Immersion_} for details).


\section{Line bundles on principal elliptic fibrations}


Let $T$ be an elliptic curve and 
$M \stackrel \pi \arrow X$ a positive principal
$T$-fibration, equipped with a preferred Hermitian
metric. Any line bundle $L$ on $M$ is stable. From
\ref{_equiva_sta_Proposition_} we obtain a natural
$\tilde T$-equivariant structure on $L$, and from 
\ref{_curva_YM_vanishes_verti_Theorem_} a connection $\nabla$
which is flat on the fibers of $\pi$. Monodromy of this
connection is related to the $\tilde T$-equivariant structure
as follows. Let $\Gamma$ be the kernel of the natural
projection $\tilde T \arrow T$, $\Gamma =\pi_1 T\cong\Z^2$.
Using the  $\tilde T$-equivariant structure, 
we find that $\Gamma$ acts on $L$ by automorphisms.
This action coincides with the monodromy of $\nabla$
along the fibers of $\pi$. 

The group of holomorphic automorphisms of $L$ is identified with
$\C^*$, therefore the monodromy map acts as a character
$\chi:\; \Gamma \arrow \C^*$. Since $\Gamma$ is 
a monodromy of a Hermitian connection, $\chi$ 
takes values in $U(1)$. We denote this character
by $\chi_L:\; \Gamma \arrow U(1)$.

\hfill

It turns out that any character can be realized by some
line bundle. The following proposition is the main result
of this section.

\hfill

\proposition \label{_line_bundles_any_char_Proposition_}
 Let $T$ be an elliptic curve and 
$M \stackrel \pi \arrow X$ a positive principal
$T$-fibration, $\dim M\geq 3$. Consider any character
$\chi:\; \Gamma \arrow U(1)$.
Then there exists a holomorphic line bundle
$L$ on $M$ such that the corresponding
character $\chi_L:\; \Gamma \arrow U(1)$
is equal to $\chi$. 

\hfill

{\bf Proof:} Consider the 
commutative diagram with exact rows coming
from the exponential exact sequence
\begin{equation}\label{_Pic_0_CD_Equation_}
\begin{CD} 
H^1(M,\calo_M) @>>> Pic_0(M) @>>> 0 \\
@VVV @VVV\\
H^1(T,\calo_T) @>>> Pic_0(T) @>>> 0 \\
\end{CD}
\end{equation}
The characters $\chi:\; \Gamma \arrow U(1)$
are in bijective correspondence with the bundles
from $Pic_0(T)$, and the correspondence is provided by
a unique flat connection on every $L\in Pic_0(T)$
(this follows e.g. from the Kobayashi-Hitchin correspondence
on elliptic curce). Therefore, to prove 
\ref{_line_bundles_any_char_Proposition_}
it is sufficient and necessary to show that 
the natural arrow $Pic_0(M) \arrow Pic_0(T)$
is surjective. As one can see from 
\eqref{_Pic_0_CD_Equation_}, this is implied
by surjectivity of the natural restriction map
\begin{equation}\label{_restri_on_H^1_O_Equation_}
H^1(M,\calo_M)\arrow H^1(T,\calo_T).
\end{equation}
 Since $\dim H^1(T,\calo_T)=1$, it is actually sufficient
to show that \eqref{_restri_on_H^1_O_Equation_}
is non-trivial. Looking again on 
\eqref{_Pic_0_CD_Equation_}, we find that
to show non-triviality of \eqref{_restri_on_H^1_O_Equation_}
it is sufficient to prove that 
$Pic_0(M) \arrow Pic_0(T)$ is non-trivial.

We reduced \ref{_line_bundles_any_char_Proposition_}
to the following lemma

\hfill

\lemma\label{_non-trivial_line_exists_Lemma_}
In assumptions of \ref{_line_bundles_any_char_Proposition_},
there exists a holomorphic line bundle $L$ on $M$, $L\in Pic_0(M)$ 
such that $L$ is non-trivial on fibers of $\pi:\; M \arrow X$.

\hfill

\remark 
Restriction of $L$ to different fibers of $\pi$ is
a flat bundle with the monodromy determined by the
character $\chi_L$. Therefore, the restriction
$L\restrict{\pi^{-1}(x)}$ is independent from 
the choice of $x\in X$ (the fibers of $\pi$ are naturally
identified because the $\pi$ is a principal fibration).

\hfill

{\bf Proof of \ref{_non-trivial_line_exists_Lemma_}:}
Consider a trivial line bundle
$L_{triv}$ on $M$ with trivial flat connection $\nabla_{triv}$.
Let $\theta$ be a 1-form which satisfies 
$d\theta=\omega_0$, where
$\omega_0 = \pi^* \omega_X$ is the pullback of
the K\"ahler form on $X$. 

Consider the connection  
$\nabla_0= \nabla_{triv}+ \theta$
on $L_{triv}$. By definition, the curvature of
$\nabla_0$ is equal $\omega_0$. Therefore,
$(\nabla_0^{0,1})^2=0$, and $\nabla_0$
defines a holomorphic structure on 
$L_{triv}$. Denote the corresponding
holomorphic line bundle by $L_0$.
The degree of $L_0$ is easy to compute:
\begin{equation} \label{_degree_line_Equation_}
  \deg L_0 = \int_M \omega \wedge \omega_0^{n-1}
\end{equation}
and this number is clearly positive. Therefore,
$L_0$ is a non-trivial holomorphic bundle. 
\ref{_non-trivial_line_exists_Lemma_} is proven.
We proved \ref{_line_bundles_any_char_Proposition_}.
\endproof

\hfill

\remark\label{_degree_triv_Remark_}
Let $q\in \R$ be a number. Consider the connection
$\nabla'= \nabla_{triv}+ q\theta$
on $L_{triv}$. Clearly, its curvature is equal
$q\omega_0$. Therefore, $\nabla'$ defines
a holomorphic structure on $L_{triv}$.
Denote the corresponding holomorphic bundle
by $L(q)$. From \eqref{_degree_line_Equation_} we
obtain that $\deg L(q) = q \deg L$.
Therefore, the trivial $C^\infty$-bundle
$L_{triv}$ admits holomorphic structures 
of any given degree.


\section{Structure theorem for stable bundles}


\subsection{Equivariant $\tilde T$-action and the stable bundles}

The main result of this paper is the following theorem.

\hfill

\theorem\label{_stru_sta_Theorem_}
Let $M\stackrel\pi\arrow X$, $\dim_\C M \geq 3$ 
be a positive principal
elliptic fibration equipped with a preferred Hermitian
metric, and $B$ a stable holomorphic bundle
on $M$. Then $B\cong L\otimes \pi^* B_0$, where
$L$ is a line bundle on $M$ and $B_0$ a stable
bundle on $X$.

\hfill

{\bf Proof:}
As \ref{_equiva_sta_Proposition_}
implies, $B$ is $\tilde T$-equivariant.
The kernel $\Gamma$ of the natural projection
$\tilde T \arrow T$ acts on $B$ by holomorphic
automorphisms. Since $B$ is stable, all its
automorphisms are proportional to identity.
Therefore, $\Gamma$ acts on $B$ by characters
$\chi:\; \Gamma \arrow \C^*$. Since 
$\chi$ can be obtained via monodromy
of the Hermitian-Einstein connection,
$\chi$ takes values in $U(1)$. Let
$L$ be a line bundle on $M$ which
is $\tilde T$-equivariant and has the
same monodromy action (such line bundle
exists by \ref{_line_bundles_any_char_Proposition_}).
Then $B\otimes L^{-1}$  has trivial monodromy on
the fibers of $\pi$. Clearly, $B\otimes L^{-1}= \pi^* B_0$
for some holomorphic bundle $B_0$ on $X$.
Since $L$ is by construction Hermitian-Einstein,
the same is true for $B\otimes L^{-1}$ and for
$B_0$ (the last assertion is true because
$\pi$ is a holomorphic Riemannian submersion).
We proved \ref{_stru_sta_Theorem_}.
\endproof

\hfill

\remark\label{_Bando_Siu_Remark_}
In \cite{_Bando_Siu_} S. Bando and Y.-T. Siu
developed Kobayashi-Hitchin correspondence for reflexive
coherent sheaves. Using these results, it is easy
to extend \ref{_stru_sta_Theorem_} for reflexive coherent sheaves.
We obtain that any stable reflexive coherent sheaf $F$
on $M$ is isomorphic to  $L\otimes \pi^* F_0$, where
$F_0$ is a stable reflexive coherent sheaf on $X$.

\subsection{Filtrable coherent sheaves on positive elliptic fibrations}

\definition
Let $Z$ be a complex variety. A coherent sheaf $F$
on $Z$ is called {\bf filtrable}
if and only if the following equivalent 
conditions are satisfied.

\begin{description}
\item[(i)] $F$ is obtained as a successive extenstion
of coherent sheaves of rank 1
\item[(ii)] $F$ admits a filtration
\[
0 = F_0 \subset F_1 \subset F_2 ... \subset F_N = F,
\]
with $\rk(F_i/F_{i-1})=1$.
\end{description}
On a quasiprojective variety $Z$, every coherent sheaf
is filtrable. This is not true if $Z$ is not algebraic,
as the following example shows.

\hfill

\example\label{_non_filtra_Example_}
Let $M$ be a K3 surface with $Pic(M)=0$ 
(generic non-algebraic K3 surfaces have 
$Pic(M)=0$). Since all line bundles on $M$
are trivial, no stable vector bundle $B$, $\rk B>1$ 
can be filtrable. However, the tangent bundle
$TM$ is stable, because it is Hermitian-Einstein.
Therefore, $TM$ is not filtrable.

\hfill

\theorem\label{_filtrable_Theorem_}
Let $M\stackrel \pi \arrow X$ be a positive principal
elliptic fibration, $\dim M\geq 3$. Assume that
$X$ is projective. Then all coherent sheaves on $M$
are filtrable.

\hfill

{\bf Proof:} Let $F$ be a coherent sheaf on $M$.
Using the Harder-Narasimhan filtration
(\cite{_Bruasse:Harder_Nara_}),
we reduce \ref{_filtrable_Theorem_}
to the case when $F$ is semistable. Using 
the Jordan-H\"older filtration,
we reduce it to the case when 
$F$ is a stable reflexive sheaf. 
By \ref{_stru_sta_Theorem_}
(see \ref{_Bando_Siu_Remark_}),
$F= \pi^* F_0 \otimes L$, where
$F_0$ is a stable reflexive 
sheaf on $X$. Since $X$ is projective,
$F_0$ is filtrable. Then $F$ is also
filtrable. We proved \ref{_filtrable_Theorem_}.
\endproof

\hfill

\hfill

{\bf Acknowledgements:} I am grateful to Ruxandra Moraru
who explained to me the geometry of holomorphic bundles on 
elliptic fibrations and gave me reference I needed, and
to Tony Pantev, who told me about
the applications and physical background of this geometry. 
Tony also found an error in a preliminary version of this paper.
My gratitude also to F. Bogomolov and D. Kaledin, 
who gave me important insights on the geometry of 
Bott-Samelson manifolds.

{\small

\hfill

\noindent {\sc Misha Verbitsky\\
University of Glasgow, Department of Mathematics, \\
15 University Gardens, Glasgow G12 8QW, Scotland.}\\
\ \\
{\sc  Institute of Theoretical and
Experimental Physics \\
B. Cheremushkinskaya, 25, Moscow, 117259, Russia }\\
\ \\
\tt verbit@maths.gla.ac.uk, \ \  verbit@mccme.ru 
}

\end{document}